\input ./style/arxiv-general.cfg
\documentclass[aos,MSNbibl,nameyear,dvips]{arximspdf}
\makeatletter
   \@ifpackageloaded{graphicx}{}{\usepackage{graphicx}}
\makeatother

\doi{10.1214/15-AOS1342}
\volume{43}
\issue{5}
\pubyear{2015}
\firstpage{2086}
\lastpage{2101}
\docsubty{FLA}

\makeatletter
\newcommand{\rrvert}{\vert}
\newcommand{\llvert}{\vert}
\newtheorem{theorem}{Theorem}[section]
\newtheorem{lemma}{Lemma}[section]
\newtheorem{corollary}{Corollary}[section]
\newproclaim{assumption}{Assumption}[section]
\newproclaim{example}{Example}[section]
\makeatother

\begin{document}
\begin{frontmatter}

\title{Adaptive testing on a regression function at~a~point}
\runtitle{Adaptive testing at a point}

\begin{aug}
\author{\fnms{Timothy}~\snm{Armstrong}\corref{}\ead[label=e1]{timothy.armstrong@yale.edu}}
\runauthor{T. Armstrong}

\affiliation{Yale University}
\address{Department of Economics\\
Yale University\\
30 Hillhouse Ave\\
New Haven, Connecticut 06511\\
USA\\
\printead{e1}}
\end{aug}

%
\received{\smonth{10} \syear{2014}}
%
\revised{\smonth{2} \syear{2015}}

%
\begin{abstract}
We consider the problem of inference on a regression function at a
point when the entire function satisfies a sign or shape restriction
under the null. We propose a test that achieves the optimal minimax
rate adaptively over a range of H\"{o}lder classes, up to a $\log\log
n$ term, which we show to be necessary for adaptation. We apply the
results to adaptive one-sided tests for the regression discontinuity
parameter under a monotonicity restriction, the value of a monotone
regression function at the boundary and the proportion of true null
hypotheses in a multiple testing problem.
\end{abstract}

%

%
\end{frontmatter}

\section{Introduction}\label{introduction_sec}

We consider a Gaussian regression model with random design. We observe
$\{(X_i,Y_i)\}_{i=1}^n$ where $X_i$ and $Y_i$ are real valued random
variables with $(X_i,Y_i)$ i.i.d. and
%
\begin{equation}
\label{reg_model} Y_i=g(X_i)+\varepsilon_i,\qquad
\varepsilon_i|X_i\sim N\bigl(0,\sigma
^2(X_i)\bigr),\qquad X_i\sim F_X,
\end{equation}
where $F_X$ denotes the c.d.f. of $X_i$.
We are interested in hypothesis tests about the regression function $g$
at a point, which we normalize to be zero.
We impose regularity conditions on
the conditional variance of $Y_i$ and the distribution of $X_i$ near
this point:
for some $\eta>0$,
%
\begin{equation}\quad\hspace*{4pt}
\label{reg_model_bdd_fs} \eta t\le\bigl\llvert F_X(t)-F_X(-t)\bigr
\rrvert \le t/\eta,\qquad \eta\le\sigma^2(x)\le1/\eta \mbox{ for $|x|<\eta$,
$0<t<\eta$.}
\end{equation}
Note that this allows (but does not impose) that our point of interest,
$0$, may be on the boundary of the support of $X_i$.

We consider the null hypotheses
%
\begin{eqnarray}
&&H_0:\bigl\{g|g(x)=0 \mbox{ all $x\in\mbox{supp}(X_1)$}
\bigr\}, \label{eq_hyp}
\\
&&H_0:\bigl\{g|g(x)\le0 \mbox{ all $x\in\mbox{supp}(X_1)$}
\bigr\}, \label
{ineq_hyp}%
\end{eqnarray}
where $\mbox{supp}(X_1)$ denotes the support of the distribution $F_X$,
and the alternative
$H_1:\{g|g(0)\ge b, g\in\mathcal{F}\}$,
where $\mathcal{F}$ imposes smoothness conditions on $g$.
In particular, we consider H\"{o}lder classes of
functions with exponent $\beta\le1$:
\[
\mathcal{F}=\Sigma(\beta,L)\equiv\bigl\{g|\bigl|g(x)-g\bigl(x'\bigr)\bigr|\le
L\bigl|x-x'\bigr|^\beta \mbox{ all $x$, $x'$}\bigr\},
\]
where $L>0$ and $0\le\beta\le1$, so that the alternative is given by
\[
H_1:g\in\mathcal{G}(b,L,\beta)\equiv\bigl\{g|g(0)\ge b
\mbox{ and } g\in \Sigma(L,\beta)\bigr\}.
\]
The focus on $g(0)$ is a normalization in the sense that the results
apply to inference on $g(x_0)$ for any point $x_0$ by redefining $X_i$
to be $X_i-x_0$, so long as the point of interest $x_0$ is known.
We also consider cases where certain shape restrictions are imposed
under the null and alternative.

For simplicity, we treat the distribution $F_X$ of $X_i$ and the
conditional variance function $\sigma^2$ as fixed and known under the
null and alternative. Thus we index probability statements with the
function $g$, which determines the joint distribution of $\{(X_i,Y_i)\}
_{i=1}^n$. We note, however, that the tests considered here can be
extended to achieve the same rates without knowledge of these
functions, so long as an upper bound for $\sup_x \sigma^2(x)$ is known
or can be estimated.

It is known
[see \citet{lepskiasymptotically2000}]
that the optimal rate for testing the null hypothesis (\ref{eq_hyp}) or
(\ref{ineq_hyp}) against the alternative $H_1$ when $g$ is known to be
in the H\"{o}lder class $\Sigma(L,\beta)$ is $n^{-\beta/(2\beta+1)}$.
That is, for any $\varepsilon>0$, there exists a constant $C_*$ such that,
for any $\alpha\in(0,1)$ and sequence of tests
$\phi_n$
with level $\alpha$ under the null hypothesis (\ref{eq_hyp}),
\[
\limsup_n \inf_{g\in\mathcal{G}(C_*n^{-\beta/(2\beta+1)},L,\beta)} E_g
\phi_n \le\alpha+\varepsilon.
\]
Furthermore, using knowledge of $\beta$, one can construct a sequence
of tests $\phi_n^*$ that are level $\alpha$ for the null hypothesis
(\ref{ineq_hyp}) [and, therefore, also level $\alpha$ for the null
hypothesis (\ref{eq_hyp})] such that,
for any $\varepsilon>0$, there exists a $C^*$ such that
%
\begin{equation}
\label{rate_lower_eq} \liminf_n \inf_{g\in\mathcal{G}(C^*n^{-\beta/(2\beta+1)},L,\beta)}
E_g \phi_n^*\ge1-\varepsilon.
\end{equation}

We ask whether a single test $\phi_n$ can achieve the rate in (\ref
{rate_lower_eq}) simultaneously for all $\beta\le1$. Such a test would
be called adaptive with respect to $\beta$. We find that the answer is
no, but that adaptivity can be obtained when the rate is modified by a
$\log\log n$ term, which we show is the necessary rate for adaptation.
In particular, we show that for $C_*$ small enough, any sequence $\phi
_n$ of level $\alpha$ tests of (\ref{eq_hyp}) must have asymptotically
trivial power for some $\beta$ in the class $\mathcal{G}(C_*(n/\log\log
n)^{-\beta/(2\beta+1)},L,\beta)$ in the sense that for any $\underline
\beta<\overline\beta\le1$,
\[
\limsup_n \inf_{\beta\in[\underline\beta,\overline\beta]}\inf
_{\mathcal
{G}(C_*(n/\log\log n)^{-\beta/(2\beta+1)},L,\beta)} E_g\phi_n\le\alpha.
\]
Furthermore, we exhibit a sequence of tests $\phi^*_n$ that achieve
asymptotic power $1$ adaptively over the classes $\mathcal{G}(C^*(n/\log
\log n)^{-\beta/(2\beta+1)}),L,\beta)$
for $C^*$ large enough,
while being level $\alpha$ for the null hypothesis (\ref{ineq_hyp}):
for any $\varepsilon>0$,
\[
\lim_{n\to\infty} \inf_{\beta\in[\varepsilon,1]} \inf
_{\mathcal
{G}(C^*(n/\log\log n)^{-\beta/(2\beta+1)},L,\beta)} E_g\phi_n^*=1.
\]

Our interest in testing at a point stems from several problems in
statistics and econometrics in which a parameter is given by the value
of a regression or density function at the boundary, and where the
function can plausibly be assumed to satisfy a monotonicity
restriction. This setup includes the regression discontinuity model and
inference on parameters that are ``identified at infinity,'' both of
which have received considerable attention in the econometrics
literature; see, among others, \citet
{chamberlainasymptotic1986,heckmanvarieties1990,andrewssemiparametric1998,hahnidentification2001}.
In the closely related problem where $g$ is a density rather than a
regression function, our setup covers the problem of inference on the
proportion of null hypotheses when testing many hypotheses; see \citet
{storeydirect2002}.
We discuss these applications in Section~\ref{applications_sec}. The
results in this paper can be used to obtain adaptive one-sided
confidence intervals for these parameters, and to show that they
achieve the minimax adaptive rate.

Problems closely related to those considered here have been considered
in the literature on asymptotic minimax bounds in nonparametric
testing, and our results draw heavily from this literature. Here, we
name only a few, and refer to \citet{ingsternonparametric2003}, for a
more thorough exposition of the literature. Typically, the goal in this
literature is to derive bounds in problems similar to the one
considered here, but with the alternative given by $\{\varphi(g)\ge b\}
\cap\mathcal{F}$, where $\varphi(g)$ is some function measuring
distance from the null and $\mathcal{F}$ a class of functions imposing
smoothness on $g$. Our problem corresponds to the case where $\varphi
(g)=g(0)$ and $\mathcal{F}=\Sigma(L,\beta)$, where we focus on
adaptivity with respect to $\beta\le1$. \citet
{lepskiasymptotically2000} consider this problem for fixed $(L,\beta
)$, and also consider the case where $\varphi(g)$ is the $\ell_\infty$ norm.
\citet{pouettesting2000} considers $\varphi(g)=g(0)$ with $\mathcal
{F}$ given by a class of analytic functions satisfying certain restrictions.
\citet{dumbgenmultiscale2001} consider the $\ell_\infty$ norm and
adaptivity over H\"{o}lder classes with respect to $(L,\beta)$ and
find, in contrast to our case, that adaptivity can be achieved without
a loss in the minimax rate (or, for adaptivity over $L$, even the
constant). In these papers, the optimal constants $C^*$ and $C_*$ are
also derived in some cases. \citet{spokoinyadaptive1996} considers
adaptivity to Besov classes under the $\ell_2$ norm and shows that, as
we derive in our case, the minimax rate can be obtained adaptively only
up to an additional $\log\log n$ term. It should also be noted that the
tests we use to achieve the minimax adaptive rate bear a close
resemblance to tests used in other adaptive testing problems; see, for
example, \citet{fantest1996,donohohigher2004}, as well as some of
the papers cited above.

Our results can be used to obtain one-sided confidence intervals for a
monotone function at the boundary of its support, which complements
results in the literature on adaptive confidence intervals for shape
restricted densities. \citet{lownonparametric1997} shows that
adaptive confidence intervals cannot be obtained without shape
restrictions on the function. \citet{caiadaptation2004} develop a
general theory of adaptive confidence intervals under shape
restrictions. \citet{caiadaptive2013} consider adaptive confidence
intervals for points on the interior of the support of a shape
restricted density and show that, in contrast to our case, the adaptive
rate can be achieved with no additional $\log\log n$ term. \citet
{dumbgenoptimal2003} considers the related problem of adaptive
confidence bands for the entire function. Our interest in points on the
boundary stems from the specific applications considered in Section~\ref{applications_sec}.

\section{Results}\label{results_sec}

We first state the lower bound for minimax adaptation.
All proofs are in Section~\ref{proofs_sec}.
For the purposes of some of the applications, we prove a slightly
stronger result in which $g$ may be known to be nonincreasing in $|x|$.
Let $\mathcal{G}_{|x|\downarrow}$ be the class of functions that are
nondecreasing on $(-\infty,0]$ and nonincreasing on $[0,\infty)$.

\begin{theorem}\label{upperbound_thm}
Let $0<\underline\beta<\overline\beta\le1$ be given.
There exists a constant $C_*$ depending only on $\underline\beta
,\overline\beta$, $L$ and the bounds on $F_X$ and $\sigma$ such that
the following holds:
Let $\phi_n$ be any sequence of tests taking the data $\{(X_i,Y_i)\}
_{i=1}^n$ to a rejection probability in $[0,1]$ with asymptotic level
$\alpha$ for the null hypothesis (\ref{eq_hyp}),
$\limsup_{n} E_0 \phi_n\le\alpha$. Then
\[
\limsup_n \inf_{\beta\in[\underline\beta,\overline\beta]}\inf
_{\mathcal
{G}(C_*(n/\log\log n)^{-\beta/(2\beta+1)},L,\beta)\cap\mathcal
{G}_{|x|\downarrow}} E_g\phi_n\le\alpha.
\]
\end{theorem}
Note that the results of the theorem imply the same results when
the requirement that $g\in\mathcal{G}_{|x|\downarrow}$ is removed from
the alternative, or when
the null is replaced by (\ref{ineq_hyp}) with the possible requirement
$g\in\mathcal{G}_{|x|\downarrow}$.

We now construct a test that achieves the $(n/\log\log n)^{\beta/(2\beta
+1)}$ rate.
For $k\in\{1,\ldots,n\}$, let $\hat g_k$ be the $k$-nearest neighbor
estimator of $g(0)$, given by
%
\begin{equation}
\label{knn_eq} \hat g_k = \frac{1}{k} \sum
_{|X_j|\le|X_{(k)}|} Y_j\qquad \mbox{where $|X_{(k)}|$ is the
$k$th least value of $|X_i|$}
\end{equation}
for $|X_{(k)}|<\eta$, and $\hat g_k=0$ otherwise, where $\eta$ is given
in (\ref{reg_model_bdd_fs}).
Let
\[
T_n=\max_{1\le k\le n} \sqrt{k}\hat g_k,
\]
and let $c_{\alpha,n}$ be the $1-\alpha$ quantile of $T_n$ under
$g(x)=0$ all $x$. Note that by the law of the iterated logarithm
[applied to the $N(0,1)$ variables $Y_i/\sigma(X_i)$ conditional on the
$X_i$'s], $\limsup_n c_{\alpha,n}/\sqrt{\log\log n}\le\sqrt{2}\sup_x
\sigma(x)$. Let $\phi_n^*$ be the test that rejects when $T_n>c_{\alpha,n}$.
%

\begin{theorem}\label{lowerbound_thm}
The test $\phi_n^*$ given above has level $\alpha$ for the null
hypothesis~(\ref{ineq_hyp}).
Furthermore, there exists a constant $C^*$ such that, for all
$\varepsilon>0$,
\[
\lim_{n\to\infty} \inf_{\beta\in[\varepsilon,1]}\inf
_{\mathcal
{G}(C^*(n/\log\log n)^{-\beta/(2\beta+1)},L,\beta)} E_g\phi_n^*%
=1.
\]
\end{theorem}

From the proof of Theorem~\ref{lowerbound_thm}, it can be seen that we
can take
$C^*=\sup_{\beta\in(0,1]} (\sqrt{2}\sup_x\sigma(x)2K^{3/2}
)^{2/(2+1/\beta)}
2^{2/(2\beta+1)} L^{1/(2\beta+1)}$
where $2/K$ is a\break lower bound on $(F_X(t)-F_X(-t))/t$.
However, this does not answer the question of the best possible
constant for the test $\phi^*_n$, or whether another test could achieve
a better constant.
While we leave these questions for future research, we briefly discuss
some conjectures.
We conjecture that, under additional regularity conditions on the
conditional variance $\sigma(x)$ and distribution of the covariate
$X_i$, a sharp constant $C(\beta,L)$ exists such
that, for arbitrary $\delta>0$, Theorem~\ref{upperbound_thm} holds with
$C_*$ replaced by $(1-\delta)C(\beta,L)$, and Theorem~\ref
{lowerbound_thm} holds with $C^*$ replaced by $(1+\delta)C(\beta,L)$
and $\phi_n^*$ replaced by a different test.
A reasonable candidate for a test statistic to achieve the optimal
constant would be a supremum over $\beta$ of normalized estimates based
on the optimal kernel given in Example~1 of \citet
{lepskiasymptotically2000}, with the bandwidth calibrated
appropriately for each $\beta$.
The conjectured behavior where minimax adaptive power goes to $\alpha$
or one on either side of a constant, where the constant does not depend
on the size $\alpha$ of the test, would be an instance of asymptotic
degeneracy related to the phenomenon observed for the $\ell_\infty$
case by \citet{lepskiasymptotically2000} (in the nonadaptive setting)
and \citet{dumbgenmultiscale2001} (for adaptivity with respect to $L$),
and our conjecture is based partly on the fact that the tests and
approximately least favorable distributions over alternatives used in
our results have a similar structure to those used in the above papers.
%

\section{Applications and extensions}\label{applications_sec}

\subsection{Inference on a monotone function at the boundary}

We note that in the case where $0$ is on the boundary of the support of
$X_i$, the results in the previous section give the optimal rate for a
one-sided test concerning $g(0)$ under a monotonicity restriction on
$g$. This can be used to obtain adaptive (up to a $\log\log n$ term)
one-sided confidence intervals for a regression function at the
boundary, where the $\log\log n$ term is necessary for adaptation. This
can be contrasted to the construction of adaptive confidence regions
for a monotone function on the interior of its support, in which case
the $\log\log n$ term is not needed; cf. \citet{caiadaptive2013}.

To form a confidence interval based on our test, we define
$T_n(\theta_0)= \max_{1\le k\le n}\sqrt{k} (\hat g_k-\theta_0
)$, and form our confidence interval by inverting tests of $H_0:g(0)\le
\theta_0$ based on $T_n(\theta_0)$ with critical value $c_{\alpha,n}$
given above (the $1-\alpha$ quantile under $g=0$ and $\theta_0=0$). The
confidence interval is then given by $[\hat c^*,\infty)$ where
$\hat c^*
=\max_{1\le k\le\overline k}[\hat g_k- c_{\alpha,n}/\sqrt{k}]$,
with $\overline k$ the largest value of $k$ such that $|X_{(k)}|<\eta$.
The following corollary to Theorems \ref{upperbound_thm} and \ref
{lowerbound_thm} shows that this CI achieves the adaptive rate.

%

%
\begin{corollary}\label{ci_thm}
Let $0<\underline\beta<\overline\beta\le1$ be given.
There exists a constant $C_*$ depending only on $\underline\beta
,\overline\beta$, $L$ and the bounds on $F_X$ and $\sigma$ such that
the following holds.
Let $[\hat c,\infty)$ be any sequence of one-sided CIs with asymptotic
coverage $1-\alpha$ for $g(0)$ when $g\in\mathcal{G}_{|x|\downarrow}$:
$\liminf_{n} \inf_{g\in\mathcal{G}_{|x|\downarrow}}P_g(g(0)\in[\hat
c,\infty))\ge1-\alpha$. Then
\[
\limsup_n \inf_{\beta\in[\underline\beta,\overline\beta]} \inf
_{g\in\Sigma(\beta,L)\cap\mathcal{G}_{|x|\downarrow}} P_g \bigl( \hat c> g(0)-C_*(n/\log\log
n)^{-\beta/(2\beta+1)} %
 \bigr)\le\alpha. %
\]
Furthermore, the CI $[\hat c^*,\infty)$ given above has coverage of at
least $1-\alpha$ for $g\in\mathcal{G}_{|x|\downarrow}$, and there
exists a $C^*$ such that, for all $\varepsilon>0$,
\[
\lim_{n\to\infty} \inf_{\beta\in[\varepsilon,1]} \inf
_{g\in\Sigma(\beta,L)\cap\mathcal{G}_{|x|\downarrow}} P_g \bigl( \hat c^*> g(0)-C^*(n/\log\log
n)^{-\beta/(2\beta+1)} %
 \bigr) =1. %
\]
\end{corollary}

The problem of inference on a regression function at the boundary has
received considerable attention in the econometrics literature, where
the problem is often termed \textit{identification at infinity}; see,
among others, \citet
{chamberlainasymptotic1986,heckmanvarieties1990,andrewssemiparametric1998,khanirregular2010}.
In such cases, it may not be plausible to assume that the density of
$X_i$ is bounded away from zero or infinity near its boundary, and the
boundary may not be finite [in which case we are interested in, e.g.,
$\lim_{x\to-\infty}g(x)$]. Such cases require one to relax the
conditions on $F_X$ in (\ref{reg_model_bdd_fs}), which can be done by
placing conditions on the behavior of $u\mapsto g(F_X^{-1}(u))$. In the
interest of space, however, we do not pursue this extension.

%

%

%

\subsection{Regression discontinuity}

Consider the regression discontinuity\break model
\[
Y_i=m(X_i)+\tau I(X_i>0)+
\varepsilon_i, \qquad\varepsilon _i|X_i\sim N
\bigl(0,\sigma^2(X_i)\bigr),\qquad X_i\sim
F_X.
\]
Here,
we strengthen (\ref{reg_model_bdd_fs}) by requiring that there exists
some $\eta>0$ such that, for all $|x|<\eta$ and $0<t<\eta$, the inequalities
$\eta t\le F_X(t)-F_X(0) \le t/\eta$,
$\eta t\le F_X(0)-F(-t) \le t/\eta$
and
$\eta\le\sigma^2(x)\le1/\eta$
are satisfied.
The regression discontinuity model has been used in a large number of
studies in empirical economics in the last decade, and has received
considerable attention in the econometrics literature; see \citet
{imbensregression2008} for a review of some of this literature.

We are interested in inference on the parameter $\tau$. Of course, $\tau
$ is not identified without constraints on $m(X_i)$. We impose a
monotonicity constraint on $m$ and ask whether a one-sided test for
$\tau$ can be constructed that is adaptive to the H\"{o}lder exponent
$\beta$ of the unknown class $\Sigma(L,\beta)$ containing $m$. In
particular, we fix $\tau_0$ and consider the null hypothesis
%
\begin{equation}
\label{rd_null} H_0: \mbox{$\tau\le\tau_0$ and $m$
nonincreasing}
\end{equation}
and the alternative
\begin{eqnarray*}
H_1: (m,\tau)&\in&\mathcal{G}^{\mathrm{rd}}(b,L,\beta)\\
& \equiv&\bigl
\{(m,\tau)|\mbox{$\tau\ge\tau_0+b$ and $m\in\Sigma(L,\beta)$
nonincreasing}\bigr\}.
\end{eqnarray*}
We extend the test of Section~\ref{results_sec} to a test that is level
$\alpha$ under $H_0$ and consistent against $H_1$ when $b=b_n$ is given
by a $\log\log n$ term times the fastest possible rate simultaneously
over $\beta\in[\varepsilon,1]$, and we show that the $\log\log n$ term
is necessary for adaptation.

To describe the test, let $\{(X_{i,1},Y_{i,1})\}_{i=1}^{n_1}$ be the
observations with $X_i\le0$, and let $\{(X_{i,2},Y_{i,2})\}
_{i=1}^{n_2}$ be the observations with $X_i> 0$. Let $\hat g_{1,k}$ be
the $k$-nearest neighbor estimator given in (\ref{knn_eq}) applied to
the sample with $X_i\le0$, and let $\hat g_{2,k}$ be defined
analogously for the sample with $X_i> 0$. Let
\[
T_{n}^{\mathrm{rd}}(\tau)=\max_{1\le k\le n}\sqrt{k}(\hat
g_{2,k}-\hat g_{1,k}-\tau).
\]
Let $c_{n,\alpha}^{\mathrm{rd}}$ be the $1-\alpha$ quantile of
$T_{n}^{\mathrm{rd}}(0)$ when $m(x)=0$ all $x$ and $\tau=0$. The test
$\phi_{n,\tau_0}^{\mathrm{rd}}$ rejects when $T_{n}^{\mathrm{rd}}(\tau
_0)>c_{n,\alpha}^{\mathrm{rd}}$.

The following corollary to Theorems \ref{upperbound_thm} and \ref
{lowerbound_thm}
gives the optimal rate for adaptive testing in the regression
discontinuity problem, and shows that the test $\phi_{n,\tau_0}^{\mathrm
{rd}}$ achieves it.
Let $E_{m,\tau}$ denote expectation under $(m,\tau)$.

%
\begin{corollary}\label{rd_thm}
Let $0<\underline\beta<\overline\beta\le1$ be given.
There exists a constant $C_*$ depending only on $\underline\beta
,\overline\beta$, $L$ and the bounds on $F_X$ and $\sigma$ such that
the following holds.
Let $\phi_n$ be any sequence of tests taking the data $\{(X_i,Y_i)\}
_{i=1}^n$ to a rejection probability in $[0,1]$ with asymptotic level
$\alpha$ for the null hypothesis (\ref{rd_null}):
$\limsup_{n} E_0 \phi_n\le\alpha$. Then
\[
\limsup_n \inf_{\beta\in[\underline\beta,\overline\beta]}\inf
_{(m,\tau
)\in\mathcal{G}^{\mathrm{rd}}(C_*(n/\log\log n)^{-\beta/(2\beta
+1)},L,\beta)} E_{m,\tau}\phi_n\le\alpha.
\]
Furthermore,
the test $\phi_{n,\tau_0}^{\mathrm{rd}}$ given above has level $\alpha$
for the null hypothesis (\ref{ineq_hyp}),
and
there exists a constant $C^*$ such that, for all $\varepsilon>0$,
\[
\lim_{n\to\infty} \inf_{\beta\in[\varepsilon,1]}\inf
_{(m,\tau)\in
\mathcal{G}^{\mathrm{rd}}(C^*(n/\log\log n)^{-\beta/(2\beta+1)},L,\beta)} E_{m,\tau}\phi_{n,\tau_0}^{\mathrm{rd}} =1.
\]
\end{corollary}
%

%

\subsection{Inference on the proportion of true null hypotheses}

Motivated by an application to large scale multiple testing, we now
consider a related setting in which we are interested in nonparametric
testing about a density, rather than a regression function. We observe
$p$-values $\{\hat p_i\}_{i=1}^n$ from $n$ independent experiments. The
$p$-values follow the mixture distribution
%
\begin{equation}
\label{fp_def_eq} \hat p_i\sim f_p(x)=\pi\cdot I\bigl(x
\in[0,1]\bigr)+(1-\pi)\cdot f_1(x),
\end{equation}
where $f_1$ is an unknown density on $[0,1]$ and $\pi$ is the
proportion of true null hypotheses.
We are interested in tests and confidence regions for $\pi$, following
a large literature on estimation and inference on $\pi$ in this
setting; see, among others, \citet{storeydirect2002}, \citet
{donohohigher2004}, \citet{meinshausenestimating2006}, \citet
{caiestimation2007} and additional references in \citet
{efronlarge-scale2012}.

Given observations from the density $f_p(x)$ with $f_1(x)$ completely
unspecified,
the best bounds that can be obtained for $\pi$ in the population are
$\pi\in[0,\overline\pi]$,
where
$\overline\pi=\overline\pi(f_p)\equiv\inf_{x\in[0,1]} f_p(x)$.
If the infimum is known to be taken at a particular location $x_0$, we
can test the null hypothesis that $\overline\pi\ge\pi_0$ against the
alternative $\overline\pi<\pi_0$ by testing the null
%
\begin{equation}
\label{mix_null} H_0: f_p(x)\ge\pi_0 \qquad\mbox{all $x$}
\end{equation}
against the alternative $f_p(x_0)<\pi_0$.
In other words,
we are interested in a version of the problem considered in Section~\ref{results_sec},
with the regression function $g$ replaced by a density
function $f_p$. Inverting these tests over $\pi_0$, we can obtain an
upper confidence interval for $\overline\pi$.
Note that since the null hypothesis $\overline\pi(f_p)\ge\pi_0$ is
equivalent to the statement that there exists a $\pi\ge\pi_0$ such
that $f_p$ follows model (\ref{fp_def_eq}) for some $f_1$, this can
also be considered a test of the null $\pi\ge\pi_0$,
and the CI can be considered a CI for $\pi$.

Assuming the $p$-values tend to be smaller when taken from the
alternative hypothesis, we can expect that $f_1(x)$ is minimized at
$x=1$ so that $f_p(x)$ will also be minimized at $1$. Following this
logic, \citet{storeydirect2002} proposes
a uniform kernel density estimator of $f_p(1)$, which can be considered
an estimator of $\overline\pi$ or of $\pi$ itself.
(In the latter case, the estimator provides an asymptotic upper bound,
but is not, in general, consistent.)
We now consider the related hypothesis testing problem
with the null given in (\ref{mix_null}) and
with the alternative
\[
H_1: f_p\in\mathcal{G}^{\pi_0}(b,L,\beta)\equiv
\bigl\{f|f_p(1)\le\pi_0-b \mbox{ and } f_p
\in\Sigma(L,\beta)\bigr\},
\]
which allows for an upper confidence interval for $\overline\pi$ (and
$\pi$ itself).
Under the maintained hypothesis that the infimum is taken at $1$,
the rate at which $b=b_n$ can approach $0$ with $H_1$ and $H_0$ being
distinguished gives the minimax rate for inference on
$\overline\pi$
when the density under the alternative is constrained to the H\"{o}lder
class $\Sigma(L,\beta)$.

To extend the approach of the previous sections to this model, let
$\hat\pi(\lambda)=\frac{1}{n(1-\lambda)}\sum_{i=1}^nI(\hat p_i>\lambda)$
be the estimate of $\pi$ used by \citet{storeydirect2002} for a given
tuning parameter $\lambda$. We form our test by searching over the
tuning parameter $\lambda$ after an appropriate normalization
\[
T_n(\pi_0)=%
\max_{0\le\lambda<1}
\sqrt{n(1-\lambda)}\bigl[\pi_0-\hat\pi(\lambda)\bigr],
\]
where we write $\max$ since the maximum is obtained.
We define our test $\phi_n(\pi_0)$ of
(\ref{mix_null})
to reject when $T_n(\pi_0)$ is greater than the critical value
$c_{n,\alpha}(\pi_0)$, given by the $1-\alpha$ quantile of $T_n(\pi_0)$
under the distribution $\pi_0\cdot\operatorname{unif}(0,1)+(1-\pi_0)\cdot
\delta_0$, where $\delta_0$ is a unit mass at $0$
and $\operatorname{unif}(0,1)$ denotes the uniform distribution on $(0,1)$.

We note that $T_n(\pi_0)$ is related to the test statistics used by
\citet{donohohigher2004} and \citet{meinshausenestimating2006}, and
can be considered a version of their approach that searches over the
larger, rather than smaller, $p$-values. \citet{donohohigher2004} set
$\pi_0=1$ and consider alternatives where $\pi$ is close to one and the
remaining $p$-values come from a normal location model with the mean
slightly perturbed, achieving a certain form of adaptivity with respect
to the amount of deviation of $\pi$ and the normal location under the
alternative. \citet{meinshausenestimating2006} consider estimation of
$\pi$ in related settings with $\pi$ close to one; see also \citet
{caiestimation2007} for additional results in this setting. In
contrast, $T_n(\pi_0)$ looks at the larger ordered $p$-values in order
to achieve adaptivity to the smoothness of the distribution of
$p$-values under the alternative in a setting where $\pi$ may not be
close to $1$.

We now state the result giving the adaptive rate for the test $\phi
_n(\pi_0)$.

\begin{theorem}\label{mix_thm}
The test $\phi_n(\pi_0)$ is level $\alpha$ for (\ref{mix_null}).
Furthermore, there exists a constant $C^*$ such that, for all
$\varepsilon>0$,
\[
\lim_{n\to\infty} \inf_{\beta\in[\varepsilon,1]} \inf
_{\mathcal{G}^{\pi_0}(C^*(n/\log\log n)^{-\beta/(2\beta+1)},L,\beta
)} E_{f_p}\phi_n=1.
\]
\end{theorem}

Given the close relation between nonparametric inference on densities
and conditional means [cf.
\citet{brownasymptotic1996},
\citet{nussbaumasymptotic1996}], a lower bound
for this problem analogous to the one given in Theorem~\ref
{upperbound_thm} for the regression problem seems likely. However, in
the interest of space, we do not pursue such an extension.

%

%

\section{Proofs}\label{proofs_sec}

\subsection{Proof of Theorem \texorpdfstring{\protect\ref{upperbound_thm}}{2.1}}

The following gives a bound on average power over certain alternatives,
and will be used to obtain a bound on minimax power over certain
alternatives conditional on $X_1,\ldots,X_n$. Note that the bound goes
to zero as $M\to\infty$ for $C<1$.

\begin{lemma}\label{il_bound_lemma}
Let $W_1,\ldots,W_N$ be independent under measures $P_0$ and $P_1,\ldots
,P_N$, with
$W_i\sim N(0,s_i^2)$ under $P_0$ and $W_i\sim N(m_{i,k},s_{i}^2)$ under $P_k$.
Let $\underline M$ and $\overline M$ be integers with $1\le
2^{\underline M}<2^{\overline M}\le N$, and let $M=\overline
M-\underline M+1$.
Let $\phi$ be a test statistic that takes the data to a rejection
probability in $[0,1]$. Suppose that for some $C$,
\[
|m_{i,k}/s_{i}| \le C \sqrt{\log M}/\sqrt{k}\qquad \mbox{all
$i,k$}
\]
and that $m_{i,k}=0$ for $i>k$.
Then
\begin{eqnarray*}
\frac{1}{M}\sum_{j=\underline M}^{\overline M}E_{P_{2^j}}
\phi- E_{P_0} \phi &\le& \sqrt{\frac{1}{M} \bigl(M^{C^2}-1
\bigr) +\frac{2}{M(\sqrt{2}-1)}C^2(\log M) M^{C^2/\sqrt{2}}}
\\
&\equiv& B(C,M). %
\end{eqnarray*}
\end{lemma}
\begin{pf}
We express the average power as the following sample mean of likelihood
ratios under the null, following arguments used in,
for example, \citet{lepskiasymptotically2000}:
\begin{eqnarray*}
\frac{1}{M}\sum_{j=\underline M}^{\overline M}E_{P_{2^j}}
\phi- E_{P_0} \phi& =&\frac{1}{M}\sum_{j=\underline M}^{\overline M}E_{P_0}
\frac
{dP_{2^j}}{dP_0}\phi- E_{P_0} \phi
\\
&=&E_{P_0} \Biggl\{\frac{1}{M}\sum
_{j=\underline M}^{\overline M} \Biggl[\exp \Biggl(\sum
_{i=1}^N\bigl(\mu_{i,j}Z_i-
\mu_{i,j}^2/2\bigr) \Biggr)-1 \Biggr]\phi \Biggr\},
\end{eqnarray*}
where $\mu_{i,j}=m_{i,2^j}/s_{i}$ and $Z_i\equiv W_i/s_{i}$ are
independent $N(0,1)$ under $P_0$.
By Cauchy--Schwarz, the above display is bounded by the square root of
%
\begin{eqnarray}
\label{sum_bound_eq} %
&&\frac{1}{M^2}\sum_{j=\underline M}^{\overline M}
\sum_{\ell=\underline
M}^{\overline M} E_{P_0} \Biggl[
\exp \Biggl(\sum_{i=1}^N\bigl(
\mu_{i,j}Z_i-\mu_{i,j}^2/2\bigr)
\Biggr)-1 \Biggr] \nonumber\\
&&\qquad\quad{}\times\Biggl[\exp \Biggl(\sum_{i=1}^N
\bigl(\mu_{i,\ell}Z_i-\mu_{i,\ell}^2/2
\bigr) \Biggr)-1 \Biggr]
\nonumber
\\
&&\qquad=\frac{1}{M^2}\sum_{j=\underline M}^{\overline M}\sum
_{\ell=\underline
M}^{\overline M} \Biggl[\exp \Biggl(\sum
_{i=1}^N\mu_{i,j}\mu_{i,\ell}
\Biggr)-1 \Biggr]
\\
&&\qquad\le\frac{1}{M^2}\sum_{j=\underline M}^{\overline M} \bigl[
\exp \bigl(C^2\log M \bigr)-1\bigr]\nonumber \\
&&\qquad\quad{}+\frac{2}{M^2}\sum
_{j=\underline M}^{\overline M}\sum_{\ell=\underline M}^{j-1}
\bigl[\exp \bigl(C^2(\log M)2^{-|j-\ell|/2} \bigr)-1\bigr],\nonumber
\end{eqnarray}
where the equality follows from using properties of the normal
distribution to evaluate the expectation,
and the last step follows by plugging in the bound
$C[\sqrt{\log M}/\sqrt{2^k}]I(i\le2^k)$ for $\mu_{i,k}=m_{i,2^k}/s_{2^k}$.
Using the fact that $\exp(x)-1\le x\cdot\exp(x)$, the inner sum of the
second term can be bounded by
\begin{eqnarray*}
&&\sum_{\ell=\underline M}^{j-1} C^2(\log
M)2^{-|j-\ell|/2} \exp \bigl(C^2(\log M)2^{-|j-\ell|/2} \bigr)
\\
&&\qquad\le C^2(\log M) \exp \bigl(C^2(\log M)/
\sqrt{2} \bigr) \sum_{k=1}^{\infty}
2^{-k/2} %
\\
&&\qquad=C^2(\log M) M^{C^2/\sqrt{2}}
\frac{1}{\sqrt{2}-1}.
\end{eqnarray*}
Plugging this into (\ref{sum_bound_eq}) and taking the square root
gives the claimed bound.
\end{pf}

Before proceeding, we recall a result regarding uniform convergence of
empirical c.d.f.s. %

\begin{lemma}\label{cdf_conv_lemma}
Let $Z_1,\ldots,Z_n$ be i.i.d. real valued random variables with c.d.f.
$F_Z$. Then, for any sequence $a_n$ with $a_n n\to\infty$,
\[
\sup_{F(z)\ge a_n} \biggl\llvert \frac{({1}/{n})\sum_{i=1}^nI(Z_i\le
z)-F_Z(z)}{F_Z(z)}\biggr\rrvert
\stackrel{p} {\to} 0.
\]
\end{lemma}

\begin{pf}
See \citet{wellnerlimit1978}, Theorem~0.
\end{pf}

Let $P_X$ denote the product measure on the $X_i$'s common to all
distributions in the model, and
let $A_n$ be the event that
%
\begin{equation}
\label{An_def_eq} %
\eta t/2 \le\frac{1}{n}\sum
_{i=1}^nI \bigl(|X_i|\le t \bigr) \le2 t/\eta
\qquad\mbox{for all }(\log n)/n <t<\eta.
\end{equation}
We will use the fact that $P_X(A_n)\to1$,
which follows by plugging condition (\ref{reg_model_bdd_fs}) into the
conclusion of Lemma~\ref{cdf_conv_lemma} for $Z_i=|X_i|$.

We now construct a function in $\mathcal{G}(b,L,\beta)$ for each $\beta
\in[\underline\beta,\overline\beta]$ that along with Lemma~\ref
{il_bound_lemma}, can be used to prove the theorem.
%

\begin{lemma}\label{upperbound_lemma}
For a given $L$, $\beta$, $n$ and $c$, define
\[
g_{\beta,n,c}(x)=\max\bigl\{c \bigl[(\log\log n)/n\bigr]^{\beta/(2\beta+1)}-L|x|^\beta
,0\bigr\}.
\]
Let $0<\underline\beta<\overline\beta$ be given.
For small enough $c$, we have the following. For any sequence of tests
$\phi_n$ taking the data into a $[0,1]$ rejection probability,
\[
\lim_{n\to\infty} \inf_{\beta\in[\underline\beta, \overline\beta]} [E_{g_{\beta,n,c}}
\phi_n-E_{0}\phi_n ] =0.
\]
\end{lemma}

\begin{pf}
Let
$\hat N(\beta)=\hat N(\beta,X_1,\ldots,X_n)
=\sum_{i=1}^n I(L|X_i|^\beta\le c[(\log\log n)/ \break n]^{\beta/(2\beta+1)})
=\sum_{i=1}^n I(|X_i| \le(c/L)^{1/\beta}[(\log\log n)/n]^{1/(2\beta+1)})$.
Let $\eta>0$ satisfy condition (\ref{reg_model_bdd_fs}).
Letting $N(\beta)= \eta^{-1}\cdot n \cdot[(\log\log n)/n]^{1/(2\beta+1)}$,
we have, for $(c/L)\le1$,
$\hat N(\beta)\le2N(\beta)$ for all $\beta\in[\underline\beta,\overline
\beta]$
on the event $A_n$ defined in (\ref{An_def_eq}).
Note that
$N(\beta)/\log\log n =\eta^{-1} (n/\log\log n)^{2\beta/(2\beta+1)}$ and
$g_{\beta,n,c}(x)\le c[(\log\log n)/n]^{\beta/(2\beta+1)}$ for all $x$,
so that
%
\begin{equation}
\label{g_bound} g_{\beta,n,c}(x)\le c\bigl[(\log\log n)/n\bigr]^{\beta/(2\beta+1)}
= c\eta^{-1/2} \bigl[ N(\beta)/\log\log n \bigr]^{-1/2}
\end{equation}
for all $x$.

Let $\underline M_n=\lceil\log_2 [2N(\underline\beta)] \rceil$ and
$\overline M_n=\lfloor\log_2 [2N(\overline\beta)] \rfloor$, and let
$\beta_{k,n}$ be such that $k=2N(\beta_{k,n})$ (so that $\underline\beta
\le\beta_{k,n}\le\overline\beta$ for $2^{\underline M_n}\le k\le
2^{\overline M_n}$).
Let $M_n=\overline M_n-\underline M_n-1$, and note that $M_n\ge(\log
n)/K$ for a constant $K$ that depends only on $\underline\beta$ and
$\overline\beta$. Plugging these into the bound in (\ref{g_bound})
yields the bound
%
\begin{equation}\qquad
\label{g_bound_k} \frac{g_{\beta_{k,n},n,c}(x)}{\sigma(x)} %
\le\frac{c\eta^{-1/2}\sqrt{2}k^{-1/2}[\log(K M_n)]^{1/2}}{\inf_{|x|<\eta}\sigma(x)} %
\le2c\eta^{-1}k^{-1/2}[\log M_n]^{1/2},
\end{equation}
where the last inequality holds for large enough $n$.
[The last equality uses the fact that $\inf_{|x|<\eta}\sigma(x)\ge\eta
^{1/2}$ for $\eta$, satisfying condition (\ref{reg_model_bdd_fs}).]

Since
$\hat N(\beta)\le2N(\beta)$ for all $\underline\beta\le\beta\le
\overline\beta$
on the event $A_n$,
we have, on this event,
letting $X_{(i)}$ be the observation $X_i$ corresponding to the $i$th
least value of $|X_i|$,
$|g_{\beta_{n,k},n,c}(X_{(i)})|=0$ for $i>2N(\beta_{n,k})=k$ for all
$\underline\beta\le\beta_{n,k}\le\overline\beta$.
Using this and the bound in (\ref{g_bound_k}), we can apply Lemma~\ref
{il_bound_lemma} conditional on $X_1,\ldots,X_n$ to obtain, for any
test $\phi$,
\[
\frac{1}{M_n} \sum_{j=\underline M_n}^{\overline M_n}
E_{g_{\beta_{n,2^j},n,c}}(\phi|X_1,\ldots,X_n)-E_{0}(
\phi|X_1,\ldots,X_n) \le B\bigl(2 c\eta^{-1},M_n
\bigr)
\]
on the event $A_n$ for large enough $n$.
Thus
\begin{eqnarray*}
&&\hspace*{-2pt}\lim_{n\to\infty} \inf_{\beta\in[\underline\beta, \overline\beta]} E_{g_{\beta,n,c}}
\phi_n-E_{0}\phi_n \\
&&\hspace*{-2pt}\qquad\le\lim
_{n\to\infty} \frac{1}{M_n} \sum_{j=\underline M_n}^{\overline M_n}
E_{g_{\beta_{n,2^j},n,c}} \phi_n-E_{0}\phi_n
\\
&&\hspace*{-2pt}\qquad\le\lim_{n\to\infty} E_{P_X} \frac{1}{M_n} \sum
_{j=\underline M_n}^{\overline M_n} \bigl[E_{g_{\beta_{n,2^j},n,c}}(
\phi|X_1,\ldots,X_n)-E_{0}(\phi
|X_1,\ldots,X_n) \bigr]I(A_n)
\\
&&\hspace*{-2pt}\qquad\quad{}+\bigl[1-P_X(A_n)\bigr]
\\
&&\hspace*{-2pt}\qquad\le\lim_{n\to\infty} B\bigl(2 c\eta^{-1},M_n
\bigr) +\bigl[1-P_X(A_n)\bigr].
\end{eqnarray*}
This converges to zero for small enough $c$.
\end{pf}

Theorem~\ref{upperbound_thm} now follows from Lemma~\ref
{upperbound_lemma} since
$g_{\beta,n,c}\in\mathcal{G}(c[(\log\log n)/ \break  n]^{\beta/ (2\beta
+1)},L,\beta)$.

\subsection{Proof of Theorem \texorpdfstring{\protect\ref{lowerbound_thm}}{2.2}}

For the test $\phi_n^*$, we have,
for $(b/L)^{1/\beta}< \eta$,
\begin{eqnarray*}
&&\inf_{g\in\mathcal{G}(b,L,\beta)} E_g \bigl(\phi_n^*|X_1,
\ldots,X_n\bigr)
\\
&&\qquad\ge\inf_{g\in\mathcal{G}(b,L,\beta)} P_g \biggl\{\frac{\sum_{|X_i|\le(b/L)^{1/\beta}} Y_i}{\sqrt{\sum_{i=1}^n I(|X_i|\le(b/L)^{1/\beta})}} >
c_{\alpha,n}|X_1,\ldots,X_n \biggr\}. %
\end{eqnarray*}
Under $P_g$, the random variable $\frac{\sum_{|X_i|\le(b/L)^{1/\beta}}
Y_i}{\sqrt{\sum_{i=1}^n I(|X_i|\le(b/L)^{1/\beta})}}$ in the
conditional probability statement above is, conditional on $X_1,\ldots
,X_n$, distributed as a normal variable with mean
%
\begin{equation}
\label{mean_bound_eq} %
\frac{\sum_{|X_i|\le(b/L)^{1/\beta}} g(X_i)}{\sqrt{\sum_{i=1}^n
I(|X_i|\le(b/L)^{1/\beta})}} \ge\frac{b}{2}
\frac{\sum_{i=1}^n I(|X_i|\le(b/(2L))^{1/\beta})}{\sqrt
{\sum_{i=1}^n I(|X_i|\le(b/L)^{1/\beta})}}
\end{equation}
and variance
%
\begin{equation}
\label{var_bound_eq} \frac{\sum_{|X_i|\le(b/L)^{1/\beta}} \sigma^2(X_i)}{\sum_{i=1}^n
I(|X_i|\le(b/L)^{1/\beta})} \le\sup_x
\sigma^2(x),
\end{equation}
where
the lower bound on the mean holds for $g\in\mathcal{G}(b,L,\beta)$
by noting that for $g\in\mathcal{G}(b,L,\beta)$, $g(x)\ge b-L|x|^\beta
$, so $g(x)\ge0$ for $|x|\le(b/L)^{1/\beta}$, and for $|x|\le
[b/(2L)]^{1/\beta}$, $g(x)\ge b-L|[b/(2L)]^{1/\beta}|^\beta=b/2$.
Let $K=2\eta^{-1}$.
On the event $A_n$ defined in (\ref{An_def_eq})
(which holds with probability approaching one),
for $b$ in the appropriate range, the right-hand side of (\ref
{mean_bound_eq}) is bounded from below by
\[
\frac{b}{2}\cdot\frac{({1}/{K})\cdot n\cdot(b/(2L))^{1/\beta}}{
\sqrt{K\cdot n\cdot(b/L)^{1/\beta}}} =\frac{1}{2K\sqrt{K}}
\cdot2^{-1/\beta}\cdot L^{-1/(2\beta)}\sqrt {n}b^{1+1/(2\beta)}.
\]
For $b=c (n/\log\log n)^{-\beta/(2\beta+1)}$, this is
\[
\frac{1}{2K\sqrt{K}}\cdot2^{-1/\beta}\cdot L^{-1/(2\beta)} c^{1+1/(2\beta)}
\sqrt{\log\log n},
\]
and for large enough $n$, this choice of $b$ is in the range where the
bounds in
(\ref{An_def_eq}),
(\ref{mean_bound_eq})
and (\ref{var_bound_eq})
can be applied for all $\beta\in[\varepsilon,1]$. Thus on the event
$A_n$,
we have for large enough $c$,
\begin{eqnarray*}
&&\inf_{\beta\in[\varepsilon,1]}\inf_{g\in\mathcal{G}(c(n/\log\log
n)^{\beta/(2\beta+1)},L,\beta)} E_g
\bigl(\phi_n^*|X_1,\ldots,X_n\bigr)
\\
&&\qquad\ge\inf_{\beta\in[\varepsilon,1]}1\\
&&\qquad\quad{}-\Phi \biggl(\frac{c_{\alpha,n}-
({1}/{(2K\sqrt{K})})\cdot2^{-1/\beta}\cdot L^{-1/(2\beta)}
c^{1+1/(2\beta)}\sqrt{\log\log n}}{\sup_x\sigma(x)} \biggr).
\end{eqnarray*}

By the law of the iterated logarithm applied to the i.i.d. $N(0,1)$
sequence $\{Y_i/\sigma(X_i)\}_{i=1}^n$,
we have $c_{\alpha,n}\le C\sqrt{\log\log n}$ for large enough $n$ for any
$C>\sqrt{2}\sup_x\sigma(x)$.
For
$c>\sup_{\beta\in(0,1]} (\sqrt{2}\sup_x\sigma(x)2K^{3/2} 2^{1/\beta
} L^{1/(2\beta)} )^{2\beta/(2\beta+1)}$,
it follows that the above display converges to $0$ as $n\to\infty$.
[Note that, while we have defined $K=2\eta^{-1}$ where $\eta$ is used
to bound both $\sigma(x)$ and $F_X(t)$ in (\ref{reg_model_bdd_fs}), $K$
is used in this proof only in bounding $F_X$ from below and can
therefore be taken to be any constant such that $2/K$ is a lower bound
on $(F_X(t)-F_X(-t))/t$ near zero.]
Since this bound holds on an event with probability approaching one,
the result follows.

%
\subsection{Proof of Corollary \texorpdfstring{\protect\ref{ci_thm}}{3.1}}

The first display follows by Theorem~\ref{upperbound_thm}
since $\phi_n=I(\hat c> 0)$ is level $\alpha$ for $H_0:g=0$, and the
display is bounded by
\[
\limsup_n\inf_{\beta\in[\underline\beta,\overline\beta]} \inf
_{g\in\mathcal{G}(C_*(n/\log\log n)^{-\beta/(2\beta+1)},L,\beta)} E_gI (\hat c > 0 ).
\]
For the second display, note that
for any constant $a$, the distribution of
$T_n(a)$ under $g$ is the same as the distribution of
$T_n(0)$ under the function $g-a$ that takes $t$ to $g(t)-a$. Thus
\[
P_g \bigl(\hat c^*> g(0)-b \bigr) %
=P_g
\bigl(T_n\bigl(g(0)-b\bigr)>c_{\alpha,n} \bigr) =P_{g-g(0)+b}
\bigl(T_n(0)>c_{\alpha,n}\bigr ).
\]
Since $g-g(0)+b$ is in $\mathcal{G}(b,L,\beta)$ for any $g\in\Sigma
(L,\beta)$, the result follows from Theorem~\ref{lowerbound_thm}.

%
\subsection{Proof of Corollary \texorpdfstring{\protect\ref{rd_thm}}{3.2}}

The proof of the second part of the corollary (the extension of Theorem~\ref{lowerbound_thm}) is similar to the original proof and is omitted.
To prove the first part of the corollary (the extension of Theorem~\ref
{upperbound_thm}), assume, without loss of generality, that $\tau_0=0$.
Define $\mbox{sgn}(X_i)$ to be $-1$ for $X_i\le0$ and $1$ for \mbox{$X_i> 0$}.
Note that, for any function $g\in\mathcal{G}(b,L/2,\beta)\cap
{G}_{|x|\downarrow}$, the function $m_g(x)=g(x)\cdot\mbox
{sgn}(X_i)-2g(0)I(X_i> 0)$ is in $\Sigma(L,\beta)$ and is
nonincreasing. [To verify H\"{o}lder continuity, note that, for $x,x'$
with $\mbox{sgn}(x)=\mbox{sgn}(x')$, $|m_g(x)-m_g(x')|\le|g(x)-g(x')|$
and, for $x,x'$ with
$\mbox{sgn}(x)\ne\mbox{sgn}(x')$,
$|m_g(x)-m_g(x')|=|g(x)-g(0)|+|g(x')-g(0)|\le(L/2)|x|^\beta
+(L/2)|x'|^\beta
\le L|x-x'|^\beta$, where the last step follows since $|x-x'|\ge x\vee x'$.]

Note that under $m=m_g$, $\tau=2g(0)$, the regression function is
$x\mapsto m_g(x)+2g(0)I(x_i>0)=g(x)\cdot\mbox{sgn}(X_i)$ so that $\{
Y_i\cdot\mbox{sgn}(X_i),X_i\}_{i=1}^n$ are distributed according to the
original regression model (\ref{reg_model}) with the given function
$g$. Of course, for $m(x)=0$ all $x$ and $\tau=0$, the regression
function is $0$ for all $x$. Thus, for any level $\alpha$ test $\phi_n$
of $(m,\tau)=(0,0)$, we can construct a test $\phi_n^*$ of (\ref
{eq_hyp}) in the original model (\ref{reg_model}) that has identical
power at $g$ to the power in the regression discontinuity model at
$(m_g,2g(0))$ for any $g$ with
$g\in\mathcal{G}(b,L/2,\beta)\cap\mathcal{G}_{|x|\downarrow}$ for some
$b,L$ and~$\beta$.
Since $(m_g,2g(0))\in\mathcal{G}^{\mathrm{rd}}(2b,L,\beta)$ whenever
$g\in\mathcal{G}(b,L/2,\beta)\cap\mathcal{G}_{|x|\downarrow}$ by the
argument above,
it follows that
\begin{eqnarray*}
&&\inf_{\beta\in[\underline\beta,\overline\beta]}
\inf_{(m,\tau)\in\mathcal{G}^{\mathrm{rd}}(2c(n/\log\log n)^{-\beta
/(2\beta+1)},L,\beta)} E_{(m,\tau)}\phi_n
\\
&&\qquad\le\inf_{\beta\in[\underline\beta,\overline\beta]}
\inf_{g\in\mathcal{G}(c(n/\log\log n)^{-\beta/(2\beta+1)},L/2,\beta
)\cap\mathcal{G}_{|x|\downarrow}}
E_{g}\phi_n^*,
\end{eqnarray*}
which converges to zero for $c$ small enough by Theorem~\ref{upperbound_thm}.

%

\subsection{Proof of Theorem \texorpdfstring{\protect\ref{mix_thm}}{3.1}}

We first show that the distribution used to obtain the critical value
is least favorable for this test statistic, so that the test does in
fact have level $\alpha$.

\begin{lemma}
The distribution $\underline f_{\pi_0}=\pi_0\cdot\operatorname{unif}(0,1)+(1-\pi_0)\delta_0$, where $\delta_0$ is a unit mass at $0$,
is least favorable for $T_n(\pi_0)$ under the null $\pi\ge\pi_0$,
\[
P_{f_p} \bigl(T_n(\pi_0)>c\bigr)\le
P_{\underline f_{\pi_0}} \bigl(T_n(\pi_0)>c\bigr) \qquad\mbox{for
$f_p$ defined by (\ref{fp_def_eq}) with $\pi \ge\pi_0$.}
\]
\end{lemma}

\begin{pf}
For $\hat p_1,\ldots,\hat p_n$ drawn from $f_p=\pi_0\cdot\operatorname{unif}(0,1)+(1-\pi_0)f_1$, let $q_1,\ldots,q_n$ be obtained from $\hat
p_1,\ldots,\hat p_n$ by setting all $\hat p_i$'s drawn from the
alternative $f_1$ to $0$. Then $T_n(\pi_0)$ weakly increases when
evaluated at the $q_i$'s instead of the $\hat p_i$'s, and
the distribution under $f_p$ of $T_n(\pi_0)$ evaluated with the $q_i$'s
is equal to the distribution under $\underline f_{\pi_0}$ of $T_n(\pi
_0)$ evaluated with the $\hat p_i$'s.
\end{pf}

The result now follows from similar arguments to the proof of
Theorem~\ref{lowerbound_thm} after noting that $c_{n,\alpha}(\pi_0)/\sqrt{\log
\log n}$ is bounded as $n\to\infty$; cf. \citet
{shorackempirical2009}, Chapter~16.

\section*{Acknowledgments}
The author would like to thank
Alexandre Tsybakov, Matias Cattaneo, Yuichi Kitamura
and Tony Cai for helpful discussions, and the Associate Editor and two
anonymous referees for helpful
comments.


\begin{thebibliography}{25}


\bibitem[\protect\citeauthoryear{Andrews and Schafgans}{1998}]{andrewssemiparametric1998}
\begin{barticle}[mr]
\bauthor{\bsnm{Andrews},~\bfnm{Donald~W.~K.}\binits{D.~W.~K.}} \AND
\bauthor{\bsnm{Schafgans},~\bfnm{Marcia~M.~A.}\binits{M.~M.~A.}}
(\byear{1998}).
\btitle{Semiparametric estimation of the intercept of a sample selection model}.
\bjournal{Rev. Econ. Stud.}
\bvolume{65}
\bpages{497--517}.
\bid{doi={10.1111/1467-937X.00055}, issn={0034-6527}, mr={1637898}}
\end{barticle}
%

\bptok{imsref}%
\endbibitem

\bibitem[\protect\citeauthoryear{Brown and Low}{1996}]{brownasymptotic1996}
\begin{barticle}[mr]
\bauthor{\bsnm{Brown},~\bfnm{Lawrence~D.}\binits{L.~D.}} \AND
\bauthor{\bsnm{Low},~\bfnm{Mark~G.}\binits{M.~G.}}
(\byear{1996}).
\btitle{Asymptotic equivalence of nonparametric regression and white noise}.
\bjournal{Ann. Statist.}
\bvolume{24}
\bpages{2384--2398}.
\bid{doi={10.1214/aos/1032181159}, issn={0090-5364}, mr={1425958}}
\end{barticle}
%

\bptok{imsref}%
\endbibitem

\bibitem[\protect\citeauthoryear{Cai, Jin and Low}{2007}]{caiestimation2007}
\begin{barticle}[mr]
\bauthor{\bsnm{Cai},~\bfnm{T.~Tony}\binits{T.~T.}},
\bauthor{\bsnm{Jin},~\bfnm{Jiashun}\binits{J.}} \AND
\bauthor{\bsnm{Low},~\bfnm{Mark~G.}\binits{M.~G.}}
(\byear{2007}).
\btitle{Estimation and confidence sets for sparse normal mixtures}.
\bjournal{Ann. Statist.}
\bvolume{35}
\bpages{2421--2449}.
\bid{doi={10.1214/009053607000000334}, issn={0090-5364}, mr={2382653}}
\end{barticle}
%

\bptok{imsref}%
\endbibitem

\bibitem[\protect\citeauthoryear{Cai and Low}{2004}]{caiadaptation2004}
\begin{barticle}[mr]
\bauthor{\bsnm{Cai},~\bfnm{T.~Tony}\binits{T.~T.}} \AND
\bauthor{\bsnm{Low},~\bfnm{Mark~G.}\binits{M.~G.}}
(\byear{2004}).
\btitle{An adaptation theory for nonparametric confidence intervals}.
\bjournal{Ann. Statist.}
\bvolume{32}
\bpages{1805--1840}.
\bid{doi={10.1214/009053604000000049}, issn={0090-5364}, mr={2102494}}
\end{barticle}
%

\bptok{imsref}%
\endbibitem

\bibitem[\protect\citeauthoryear{Cai, Low and Xia}{2013}]{caiadaptive2013}
\begin{barticle}[mr]
\bauthor{\bsnm{Cai},~\bfnm{T.~Tony}\binits{T.~T.}},
\bauthor{\bsnm{Low},~\bfnm{Mark~G.}\binits{M.~G.}} \AND
\bauthor{\bsnm{Xia},~\bfnm{Yin}\binits{Y.}}
(\byear{2013}).
\btitle{Adaptive confidence intervals for regression functions under shape constraints}.
\bjournal{Ann. Statist.}
\bvolume{41}
\bpages{722--750}.
\bid{doi={10.1214/12-AOS1068}, issn={0090-5364}, mr={3099119}}
\end{barticle}
%

\bptok{imsref}%
\endbibitem

\bibitem[\protect\citeauthoryear{Chamberlain}{1986}]{chamberlainasymptotic1986}
\begin{barticle}[mr]
\bauthor{\bsnm{Chamberlain},~\bfnm{Gary}\binits{G.}}
(\byear{1986}).
\btitle{Asymptotic efficiency in semiparametric models with censoring}.
\bjournal{J.~Econometrics}
\bvolume{32}
\bpages{189--218}.
\bid{doi={10.1016/0304-4076(86)90038-2}, issn={0304-4076}, mr={0864926}}
\end{barticle}
%

\bptok{imsref}%
\endbibitem

\bibitem[\protect\citeauthoryear{Donoho and Jin}{2004}]{donohohigher2004}
\begin{barticle}[mr]
\bauthor{\bsnm{Donoho},~\bfnm{David}\binits{D.}} \AND
\bauthor{\bsnm{Jin},~\bfnm{Jiashun}\binits{J.}}
(\byear{2004}).
\btitle{Higher criticism for detecting sparse heterogeneous mixtures}.
\bjournal{Ann. Statist.}
\bvolume{32}
\bpages{962--994}.
\bid{doi={10.1214/009053604000000265}, issn={0090-5364}, mr={2065195}}
\end{barticle}
%

\bptok{imsref}%
\endbibitem

\bibitem[\protect\citeauthoryear{D{\"u}mbgen}{2003}]{dumbgenoptimal2003}
\begin{barticle}[mr]
\bauthor{\bsnm{D{\"u}mbgen},~\bfnm{Lutz}\binits{L.}}
(\byear{2003}).
\btitle{Optimal confidence bands for shape-restricted curves}.
\bjournal{Bernoulli}
\bvolume{9}
\bpages{423--449}.
\bid{doi={10.3150/bj/1065444812}, issn={1350-7265}, mr={1997491}}
\end{barticle}
%

\bptok{imsref}%
\endbibitem

\bibitem[\protect\citeauthoryear{D{\"u}mbgen and Spokoiny}{2001}]{dumbgenmultiscale2001}
\begin{barticle}[mr]
\bauthor{\bsnm{D{\"u}mbgen},~\bfnm{Lutz}\binits{L.}} \AND
\bauthor{\bsnm{Spokoiny},~\bfnm{Vladimir~G.}\binits{V.~G.}}
(\byear{2001}).
\btitle{Multiscale testing of qualitative hypotheses}.
\bjournal{Ann. Statist.}
\bvolume{29}
\bpages{124--152}.
\bid{doi={10.1214/aos/996986504}, issn={0090-5364}, mr={1833961}}
\end{barticle}
%

\bptok{imsref}%
\endbibitem

\bibitem[\protect\citeauthoryear{Efron}{2010}]{efronlarge-scale2012}
\begin{bbook}[mr]
\bauthor{\bsnm{Efron},~\bfnm{Bradley}\binits{B.}}
(\byear{2010}).
\btitle{Large-Scale Inference: Empirical Bayes Methods for Estimation, Testing, and Prediction}.
\bseries{Institute of Mathematical Statistics (IMS) Monographs}
\bvolume{1}.
\bpublisher{Cambridge Univ. Press},
\blocation{Cambridge}.
\bid{doi={10.1017/CBO9780511761362}, mr={2724758}}
\bptnote{check year}%
\end{bbook}
%

\bptok{imsref}%
\endbibitem

\bibitem[\protect\citeauthoryear{Fan}{1996}]{fantest1996}
\begin{barticle}[mr]
\bauthor{\bsnm{Fan},~\bfnm{Jianqing}\binits{J.}}
(\byear{1996}).
\btitle{Test of significance based on wavelet thresholding and {N}eyman's truncation}.
\bjournal{J. Amer. Statist. Assoc.}
\bvolume{91}
\bpages{674--688}.
\bid{doi={10.2307/2291663}, issn={0162-1459}, mr={1395735}}
\end{barticle}
%

\bptok{imsref}%
\endbibitem

\bibitem[\protect\citeauthoryear{Hahn, Todd and Van~der Klaauw}{2001}]{hahnidentification2001}
\begin{barticle}[author]
\bauthor{\bsnm{Hahn},~\bfnm{Jinyong}\binits{J.}},
\bauthor{\bsnm{Todd},~\bfnm{Petra}\binits{P.}} \AND
\bauthor{\bsnm{Van~der Klaauw},~\bfnm{Wilbert}\binits{W.}}
(\byear{2001}).
\btitle{Identification and estimation of treatment {effects} with a regression-{discontinuity} {design}}.
\bjournal{Econometrica}
\bvolume{69}
\bpages{201--209}.
\end{barticle}
%

\bptok{imsref}%
\endbibitem

\bibitem[\protect\citeauthoryear{Heckman}{1990}]{heckmanvarieties1990}
\begin{barticle}[author]
\bauthor{\bsnm{Heckman},~\bfnm{James}\binits{J.}}
(\byear{1990}).
\btitle{Varieties of selection {Bias}}.
\bjournal{The American Economic Review}
\bvolume{80}
\bpages{313--318}.
\end{barticle}
%

\bptok{imsref}%
\endbibitem

\bibitem[\protect\citeauthoryear{Imbens and Lemieux}{2008}]{imbensregression2008}
\begin{barticle}[mr]
\bauthor{\bsnm{Imbens},~\bfnm{Guido~W.}\binits{G.~W.}} \AND
\bauthor{\bsnm{Lemieux},~\bfnm{Thomas}\binits{T.}}
(\byear{2008}).
\btitle{Regression discontinuity designs: A guide to practice}.
\bjournal{J.~Econometrics}
\bvolume{142}
\bpages{615--635}.
\bid{doi={10.1016/j.jeconom.2007.05.001}, issn={0304-4076}, mr={2416821}}
\end{barticle}
%

\bptok{imsref}%
\endbibitem

\bibitem[\protect\citeauthoryear{Ingster and Suslina}{2003}]{ingsternonparametric2003}
\begin{bbook}[mr]
\bauthor{\bsnm{Ingster},~\bfnm{Y.~I.}\binits{Y.~I.}} \AND
\bauthor{\bsnm{Suslina},~\bfnm{I.~A.}\binits{I.~A.}}
(\byear{2003}).
\btitle{Nonparametric Goodness-of-Fit Testing Under {G}aussian Models}.
\bseries{Lecture Notes in Statistics}
\bvolume{169}.
\bpublisher{Springer},
\blocation{New York}.
\bid{doi={10.1007/978-0-387-21580-8}, mr={1991446}}
\end{bbook}
%

\bptok{imsref}%
\endbibitem

\bibitem[\protect\citeauthoryear{Khan and Tamer}{2010}]{khanirregular2010}
\begin{barticle}[mr]
\bauthor{\bsnm{Khan},~\bfnm{Shakeeb}\binits{S.}} \AND
\bauthor{\bsnm{Tamer},~\bfnm{Elie}\binits{E.}}
(\byear{2010}).
\btitle{Irregular identification, support conditions, and inverse weight estimation}.
\bjournal{Econometrica}
\bvolume{78}
\bpages{2021--2042}.
\bid{doi={10.3982/ECTA7372}, issn={0012-9682}, mr={2768989}}
\end{barticle}
%

\bptok{imsref}%
\endbibitem

\bibitem[\protect\citeauthoryear{Lepski and Tsybakov}{2000}]{lepskiasymptotically2000}
\begin{barticle}[mr]
\bauthor{\bsnm{Lepski},~\bfnm{O.~V.}\binits{O.~V.}} \AND
\bauthor{\bsnm{Tsybakov},~\bfnm{A.~B.}\binits{A.~B.}}
(\byear{2000}).
\btitle{Asymptotically exact nonparametric hypothesis testing in sup-norm and at a fixed point}.
\bjournal{Probab. Theory Related Fields}
\bvolume{117}
\bpages{17--48}.
\bid{doi={10.1007/s004400050265}, issn={0178-8051}, mr={1759508}}
\end{barticle}
%

\bptok{imsref}%
\endbibitem

\bibitem[\protect\citeauthoryear{Low}{1997}]{lownonparametric1997}
\begin{barticle}[mr]
\bauthor{\bsnm{Low},~\bfnm{Mark~G.}\binits{M.~G.}}
(\byear{1997}).
\btitle{On nonparametric confidence intervals}.
\bjournal{Ann. Statist.}
\bvolume{25}
\bpages{2547--2554}.
\bid{doi={10.1214/aos/1030741084}, issn={0090-5364}, mr={1604412}}
\end{barticle}
%

\bptok{imsref}%
\endbibitem

\bibitem[\protect\citeauthoryear{Meinshausen and Rice}{2006}]{meinshausenestimating2006}
\begin{barticle}[mr]
\bauthor{\bsnm{Meinshausen},~\bfnm{Nicolai}\binits{N.}} \AND
\bauthor{\bsnm{Rice},~\bfnm{John}\binits{J.}}
(\byear{2006}).
\btitle{Estimating the proportion of false null hypotheses among a large number of independently tested hypotheses}.
\bjournal{Ann. Statist.}
\bvolume{34}
\bpages{373--393}.
\bid{doi={10.1214/009053605000000741}, issn={0090-5364}, mr={2275246}}
\end{barticle}
%

\bptok{imsref}%
\endbibitem

\bibitem[\protect\citeauthoryear{Nussbaum}{1996}]{nussbaumasymptotic1996}
\begin{barticle}[mr]
\bauthor{\bsnm{Nussbaum},~\bfnm{Michael}\binits{M.}}
(\byear{1996}).
\btitle{Asymptotic equivalence of density estimation and {G}aussian white noise}.
\bjournal{Ann. Statist.}
\bvolume{24}
\bpages{2399--2430}.
\bid{doi={10.1214/aos/1032181160}, issn={0090-5364}, mr={1425959}}
\end{barticle}
%

\bptok{imsref}%
\endbibitem

\bibitem[\protect\citeauthoryear{Pouet}{1999}]{pouettesting2000}
\begin{barticle}[mr]
\bauthor{\bsnm{Pouet},~\bfnm{Christophe}\binits{C.}}
(\byear{1999}).
\btitle{On testing nonparametric hypotheses for analytic regression functions in {G}aussian noise}.
\bjournal{Math. Methods Statist.}
\bvolume{8}
\bpages{536--549}.
\bid{issn={1066-5307}, mr={1755899}}
\bptnote{check pages, check year}%
\end{barticle}
%

\bptok{imsref}%
\endbibitem

\bibitem[\protect\citeauthoryear{Shorack and Wellner}{2009}]{shorackempirical2009}
\begin{bbook}[author]
\bauthor{\bsnm{Shorack},~\bfnm{Galen~R.}\binits{G.~R.}} \AND
\bauthor{\bsnm{Wellner},~\bfnm{Jon~A.}\binits{J.~A.}}
(\byear{2009}).
\btitle{Empirical Processes with Applications to {Statistics}}.
\bpublisher{SIAM},
\blocation{Philadelphia, PA}.
\end{bbook}
%

\bptok{imsref}%
\endbibitem

\bibitem[\protect\citeauthoryear{Spokoiny}{1996}]{spokoinyadaptive1996}
\begin{barticle}[mr]
\bauthor{\bsnm{Spokoiny},~\bfnm{V.~G.}\binits{V.~G.}}
(\byear{1996}).
\btitle{Adaptive hypothesis testing using wavelets}.
\bjournal{Ann. Statist.}
\bvolume{24}
\bpages{2477--2498}.
\bid{doi={10.1214/aos/1032181163}, issn={0090-5364}, mr={1425962}}
\end{barticle}
%

\bptok{imsref}%
\endbibitem

\bibitem[\protect\citeauthoryear{Storey}{2002}]{storeydirect2002}
\begin{barticle}[mr]
\bauthor{\bsnm{Storey},~\bfnm{John~D.}\binits{J.~D.}}
(\byear{2002}).
\btitle{A direct approach to false discovery rates}.
\bjournal{J. R. Stat. Soc. Ser. B. Stat. Methodol.}
\bvolume{64}
\bpages{479--498}.
\bid{doi={10.1111/1467-9868.00346}, issn={1369-7412}, mr={1924302}}
\end{barticle}
%

\bptok{imsref}%
\endbibitem

\bibitem[\protect\citeauthoryear{Wellner}{1978}]{wellnerlimit1978}
\begin{barticle}[mr]
\bauthor{\bsnm{Wellner},~\bfnm{Jon~A.}\binits{J.~A.}}
(\byear{1978}).
\btitle{Limit theorems for the ratio of the empirical distribution function to the true distribution function}.
\bjournal{Z. Wahrsch. Verw. Gebiete}
\bvolume{45}
\bpages{73--88}.
\bid{issn={0178-8051}, mr={0651392}}
\end{barticle}
%

\bptok{imsref}%
\endbibitem
\end{thebibliography}
%





\printaddresses
\end{document}